\documentclass{amsproc}
\usepackage{graphicx}
\usepackage{amssymb}
\usepackage{epstopdf}
\usepackage{amsmath,amsthm,amscd,amssymb}
\usepackage{latexsym}
\usepackage[colorlinks,citecolor=red,pagebackref,hypertexnames=false]{hyperref}
\usepackage{geometry}                
\numberwithin{equation}{section}

\theoremstyle{plain}
\newtheorem{theorem}{Theorem}[section]
\newtheorem{lemma}[theorem]{Lemma}

\theoremstyle{definition}
\newtheorem{definition}[theorem]{Definition}

\newtheorem{case[theorem]}{Case}

\def\dim{\textrm{ dim }}

\def\R{\mathbb R}

\def\vx{\vec{x}}
\def\vy{\vec{y}}
\def\hx{\hat{x}}
\def\hy{\hat{y}}

\def\V{\mathcal V}
\def\hd{dim_{\mathcal H}}
\def\hde{\hd(E)}
\def\R{\Bbb R}

\theoremstyle{remark}
\newtheorem{remark}[theorem]{Remark}

\numberwithin{equation}{section}

\begin{document}

\title{\parbox{14cm}{\centering{On volumes determined by subsets of Euclidean space}}}


\author{Allan Greenleaf, Alex Iosevich and Mihalis Mourgoglou}

\address{Department of Mathematics \\ University of Rochester\\ Rochester, NY 14627}
\email{allan@math.rochester.edu}
\address{Department of Mathematics \\ University of Rochester\\ Rochester, NY 14627}
\email{iosevich@math.rochester.edu}
\address{Department of Mathematics, Universit\'{e} Paris-Sud 11, Orsay }
\email{mihalis.mourgoglou@math.u-psud.fr}

\thanks{The first two authors authors were supported by NSF grants DMS-0853892 and DMS-1045404. The third author holds a Sophie Germain International post-doctoral scholarship at  Fondation de Math\'{e}matiques Jacques Hadamard (FMJH) and would like to thank the faculty and staff of the Universit\'{e} Paris-Sud 11, Orsay for their hospitality. }

\begin{abstract} Given $E \subset {\Bbb R}^d$, define the \emph{volume set} of $E$, ${\mathcal V}(E)= \{det(x^1, x^2, \dots x^d): x^j \in E\}$. In $\R^3$, we prove that ${\mathcal V}(E)$ has positive Lebesgue measure if either  the Hausdorff dimension of $E\subset \Bbb R^3$ is greater than $\frac{13}{5}$, or $E$ is a product set of the form $E=B_1\times B_2\times B_3$ with $B_j\subset\R,\, dim_{\mathcal H}(B_j)>\frac23,\, j=1,2,3$. We show that the same conclusion holds for  $\V(E)$ of  Salem subsets  $E\subset\R^d$ with $\hde>d-1$, and give applications  to discrete combinatorial geometry.

\end{abstract} 

\maketitle


\section{Introduction}

\vskip.125in 

A large class of Erd\H os type problems in geometric combinatorics ask whether a large set of points in Euclidean space determines a suitable large sets of geometric relations, configurations or objects. For example, the classical Erd\H os distance problem asks whether $N$ points in ${\Bbb R}^d$, $d \ge 2$, determines $\gtrapprox N^{\frac{2}{d}}$ distinct distances. 
See, for example \cite{BMP00,Ma02,P05,PS04,Sz97} and the references  there for  thorough descriptions of these types of problems and recent results. 
(Here, and throughout, $X \lessapprox Y$, if the controlling parameter is $N$, means that for every $\epsilon>0$ there exists $C_{\epsilon}>0$ such that $X \leq C_{\epsilon} N^{\epsilon} Y$, while $X\lesssim Y$ means $X\le CY$ with $C$ independent of $N$.) 

Continuous variants of Erd\H os type geometric problems have also received much attention in recent decades. Perhaps the best known of these is the Falconer distance problem \cite{Fal86}, which asks whether the one-dimensional Lebesgue measure of the distance set $\{|x-y|: x,y \in E\}$ is positive, provided that the Hausdorff dimension, $\hde$, of $E \subset {\Bbb R}^d$, $d \ge 2$, is greater than $\frac{d}{2}$. See \cite{W99,Erd05} for the best currently known results on this problem. Also see \cite{CEHIT10} for the closely related problem of finite point configurations. 

In this paper we study the sets of \emph{volumes} determined by sets $E \subset {\Bbb R}^d$. Given $d$ vectors $x^1, x^2, \dots x^d$ in ${\Bbb R}^d$, $d \ge 2$,  let $det(x^1, \dots, x^d)$ denote the determinant of the matrix whose $j$th column is $x^j$. For $E \subset {[0,1]}^d$,  define the \emph{volume set} of $E$,

$$\V(E):=\{ det(x^1, x^2, \dots, x^d)\in\R: x^j \in E\}. $$ 

A problem is to find the optimal theshhold $s_\V(d)$ such that if the Hausdorff dimension of $E \subset {[0,1]}^d$ is greater than $s_\V(d)$, then the Lebesgue measure of $\V(E)$, denoted by ${\mathcal L}^1(\V(E))$, is positive. Letting $E$  be  a $(d-1)$-dimensional hyperplane shows that one must take $s_\V(d)\ge d-1$. 

A result due to Erdo\~{g}an, Hart and the second author \cite{EHI11} shows that  $s_\V(2)\le\frac{3}{2}$. More generally, they prove the following result. 

\begin{theorem} \label{witherdogan} \cite{EHI11} For $d \ge 2$, let $E, F \subset {[0,1]}^d$  and  $Q$ be a non-degenerate bilinear form on $\R^d$. Suppose that 
$ dim_{{\mathcal H}}(E)+\dim_{{\mathcal H}}(F)>d+1$. Then ${\mathcal L}^1(\{Q(x,y): x \in E, y \in F\})>0.$  

Thus, if $E\subset\R^2$ with $\hde>\frac32$, then $\mathcal L^1(\V(E))>0$.
\end{theorem} 

Our four main results are the following. 

\begin{theorem} \label{maingeneral}  Let $E \subset {[0,1]}^3$, with $\hde>\frac{13}{5}$. Then ${\mathcal L}^1(\V(E))>0$. 
\end{theorem} 

In view of Thm. \ref{witherdogan} above, Thm. \ref{maingeneral} may be reduced to studying the Hausdorff dimension of the set of $(d-1)$-vectors of $E$,

\begin{equation} \label{wedge} \Lambda(E):=\{\, *(x^1 \wedge x^2 \wedge \dots \wedge x^{d-1}): x^j \in E,\, 1\le j\le d-1\}, \end{equation} 
where $*$ is the Hodge star operator, $*:\Lambda^{d-1}\R^d\longrightarrow \R^d$.
Since

\begin{equation} \label{generalizedcrossproduct} det(x^1,\dots,x^d)=\pm x^d\cdot * (x^1 \wedge x^2 \wedge \dots \wedge x^{d-1}), \end{equation} 
establishing Thm. \ref{maingeneral} reduces to proving the following result. 

\begin{theorem} \label{crossproduct} Let $E \subset {[0,1]}^3$ and suppose that  $\hde >\frac{13}{5}$. Then 

\begin{equation} \label{feeder} dim_{{\mathcal H}}(E)+dim_{{\mathcal H}}(\Lambda(E))>4.\end{equation}

\end{theorem}

If one could strengthen 
(\ref{feeder}) to say that  the left hand side  is greater than  $3$, then Thm. \ref{maingeneral} would be improved to the $\hde>2$. This would be optimal because, as  noted above, the conclusion of Thm. \ref{maingeneral} does not in general hold if  $\hde\le d-1$. 

\begin{remark} In the  context of three dimensional vector spaces over finite fields, the analogue of Theorem \ref{maingeneral} is completely resolved, with a sharp exponent using sum-product technology; see \cite{CHIKR10}. However, sum-product issues are generally believed to be more difficult in the continuous setting. \end{remark} 

It is possible to obtain a better exponent than in Thm. \ref{crossproduct}, and thus Thm. \ref{maingeneral}, if the  set $E$ has a special form or satisfies Fourier decay conditions. We first consider the situation where the set under consideration is a Cartesian product of subsets of the real line. 

\begin{theorem} \label{cartesianproduct} Suppose that $E=B_1 \times B_2 \times B_3$, where  $B_j \subset [0,1]$ satisfy $dim_{\mathcal H}(B_j)>\frac{2}{3},\, j=1,2,3$. Then ${\mathcal L^1}({\mathcal V}(E))>0$. 
\end{theorem} 

Note that if each $B_j$ is Ahlfors-David regular (see e.g. \cite{M95}), the assumption that the Hausdorff dimension of each $B_j$ is greater than $\frac{2}{3}$ is equivalent to the assumption that the Hausdorff dimension of $E=B_1 \times B_2 \times B_3$ is greater than $2$. See, e.g., \cite{Falc86,M95}. 

Alternatively, one can also consider the situation in higher dimensions, when  set $E\subset\R^d,\, d\ge2$, supports a measure that has good Fourier decay properties.  Recall that $E \subset {[0,1]}^d$, $d \ge 2$, is \emph{Salem} if there exists a probability measure $\mu$ supported on $E$ such that
\begin{equation} \label{salem} |\widehat{\mu}(\xi)| \lessapprox {|\xi|}^{-\frac{s}{2}}, \quad |\xi|\longrightarrow\infty,
\end{equation} 
where $s$ denotes the Hausdorff dimension of $E$ and 
$$ \widehat{\mu}(\xi)=\int e^{-2 \pi i x \cdot \xi} d\mu(x)$$ 
is the Fourier transform of $\mu$. We prove

\begin{theorem} \label{mainsalem} 

Suppose $E\subset \R^d$ is Salem and  $s:=\hde>d-1$. Then ${\mathcal L}^1(\V(E))>0$. 
\end{theorem} 

We note that while one cannot do better than the exponent $d-1$, as the example of a hyperplane shows, it is possible that this might not be the case for Salem sets. It is  conceivable that the conclusion of Theorem \ref{mainsalem} holds if the Hausdorff dimension of $E$ is merely greater than one. 

\vskip.25in 

\section{Proof of Theorem \ref{mainsalem}} 

\vskip.125in 

Let $\psi$ be a smooth cut-off function on $\R^d$, $\equiv 1$ in the unit ball and supported in the ball of radius $\sqrt{d}$. Then 

$$ (\mu \times \mu \times \dots \times \mu) \{(x^1, x^2, \dots, x^d): t \leq det(x^1, x^2, \dots , x^d) \leq t+\epsilon \}$$ 

$$ \approx \int \int \dots \int \psi\left( \frac{det(x^1, \dots, x^d)-t}{\epsilon} \right) d\mu(x^1) \dots d\mu(x^d)$$

\begin{equation} \label{initialreduction}= \int \int \dots \int \int e^{2 \pi i \frac{\lambda}{\epsilon} det(x^1, \dots, x^d)} e^{-2 \pi i \frac{t \lambda}{\epsilon}} d\mu(x^1) \dots d\mu(x^d) \widehat{\psi}(\lambda) d\lambda.\end{equation}

Using (\ref{generalizedcrossproduct}),
it follows that (\ref{initialreduction}) equals 

$$ \epsilon \int \dots \int \widehat{\mu}(\lambda(x^1 \wedge \dots \wedge x^{d-1})) e^{-2 \pi i t \lambda} d\mu(x^1) \dots d\mu(x^{d-1}) \widehat{\psi}(\epsilon \lambda) d\lambda$$
$$=\epsilon \int \dots \int_{|x^1 \wedge \dots \wedge x^{d-1}| \leq \lambda^{-1}} 
\widehat{\mu}(\lambda(x^1 \wedge \dots \wedge x^{d-1})) e^{-2 \pi i t \lambda} d\mu(x^1) \dots d\mu(x^{d-1}) \widehat{\psi}(\epsilon \lambda) d\lambda$$
$$+\epsilon \int \dots \int_{|x^1 \wedge \dots \wedge x^{d-1}|>\lambda^{-1}} \widehat{\mu}(\lambda(x^1 \wedge \dots \wedge x^{d-1})) e^{-2 \pi i t \lambda} d\mu(x^1) \dots d\mu(x^{d-1}) \widehat{\psi}(\epsilon \lambda) d\lambda:=I+II.$$

\vskip.125in 

Since $|\widehat{\mu}(\xi)| \leq 1$, we have 

$$ I \lesssim \epsilon \int (\mu \times \dots \times \mu) \{(x^1, \dots, x^{d-1}): |x^1 \wedge \dots \wedge x^{d-1}| \leq \lambda^{-1} \} |\widehat{\psi}(\epsilon \lambda)| d\lambda.$$

\begin{lemma} \label{crossproductestimation} With the notation above, 

\begin{equation} \label{volumetrivial} (\mu \times \dots \times \mu) \{(x^1, \dots, x^{d-1}): |x^1 \wedge \dots \wedge x^{d-1}| \leq \lambda^{-1} \} \lesssim \min \{1, \lambda^{d-2-s}\}. \end{equation} 

\end{lemma} 

To prove the lemma, start by noting that one has
$$ |x^1 \wedge x^2 \wedge \dots \wedge x^{d-1}|=|x^1| \cdot |x^2 \wedge \dots \wedge x^{d-1}|\cdot \sin(\theta),$$ 
where $\theta$ is the angle between $x^1$ and the $(d-2)$-plane spanned $\Pi$ by $x^2,\dots,x^{d-1}$.
Localize to where $|x^2 \wedge \dots \wedge x^{d-1}| \approx 2^{-j}$. By induction, the measure of this set is $\lesssim 2^{j(d-3-s)}$. Since $x^1$ is contained in a $2^j \lambda^{-1}$-tubular nhood of $\Pi$, it follows that 
\begin{eqnarray*} (\mu \times \dots \times \mu) \{(x^1, \dots, x^{d-1}): |x^1 \wedge \dots \wedge x^{d-1}| \leq \lambda^{-1} \} &\lesssim& {(2^j \lambda^{-1})}^{-(d-2)} {(2^j \lambda^{-1})}^s {(2^{-j})}^{s-(d-3)}\ \\
&=&\lambda^{d-2-s} 2^{-j}.
\end{eqnarray*} 
It follows that the left hand side of (\ref{volumetrivial}) is bounded by a constant multiple of 
$$\sum_j \lambda^{d-2-s} 2^{-j} \lesssim \lambda^{d-2-s},$$  
completing the proof of Lemma \ref{crossproductestimation}. 

\vskip.125in 

It follows from Lemma \ref{crossproductestimation} that 
$$ I \lesssim \int \min \{1, \lambda^{d-2-s} \} |\widehat{\psi}(\epsilon \lambda)|\, d\lambda,$$ 
which is $\lesssim 1$ if $s>d-1$.  By the Salem property (\ref{salem}), 

$$ II \lessapprox \epsilon \int \dots \int \lambda^{-\frac{s}{2}} {|x^1 \wedge \dots \wedge x^{d-1}|}^{-\frac{s}{2}} d\mu(x^1) d\mu(x^2) \dots d\mu(x^{d-1}) |\widehat{\psi}(\epsilon \lambda)| d\lambda$$

$$ \approx \sum_{j \leq \log_2(\lambda)} 2^{\frac{js}{2}} 
\epsilon \int \lambda^{-\frac{s}{2}} \mu \times \dots \times \mu \{(x^1, \dots, x^{d-1}): |x^1 \wedge \dots \wedge x^{d-1}| \approx 2^{-j} \} \ |\widehat{\psi}(\epsilon \lambda)| \ d\lambda.$$ 
Using Lemma \ref{crossproductestimation}, one sees that this is this  is 

$$ \approx \int \sum_{j \leq \log_2(\lambda)} 2^{j(d-2-s)} 2^{\frac{js}{2}} \lambda^{-\frac{s}{2}} |\widehat{\psi}(\epsilon \lambda)| d\lambda \lesssim \int_1^{\infty} \lambda^{d-2-s} |\widehat{\psi}(\epsilon \lambda)| d\lambda,$$ 
which is bounded independently of $\epsilon$ if $s>d-1$, as desired. This finishes the proof of Thm. \ref{mainsalem}.

\vskip.25in 

\section{Proof of Theorem \ref{cartesianproduct}} 

\vskip.125in 

We shall need the following consequence of Theorem 1.0.3 in \cite{EHI11}. 
\begin{theorem} \label{tube} \cite{EHI11} Let $A_1,A_2,A_3,A_4 \subset [0,1]$, each with $\hd(A_j)>\frac{2}{3}$. Then 
\begin{equation} \label{cartesianorgasm} {\mathcal L}^1(\{a_1a_2+a_3a_4: a_j \in A_j \})>0. \end{equation} 
\end{theorem} 

\vskip.125in 

To apply Thm. \ref{tube}, consider 
\begin{equation} \label{8thgrade} det(x,y,z)=y_1(x_3z_2-x_2z_3)+y_2(x_1z_3-x_3z_1)+y_3(-x_1z_2-x_2z_1).\end{equation}  

Fix $x_1,x_2,x_3, y_3, z_3\in A$, all $\ne 0$, and let 
$$w_1=x_3z_2-x_2z_3, \ w_2=x_1z_3-x_3z_1.$$ Observe that 
$$ -x_1z_2-x_2z_1=-\frac{x_1}{x_3}w_1-\frac{x_2}{x_3}w_2.$$ 

It follows that the expression in (\ref{8thgrade}) equals 
$$ y_1w_1+y_2w_2-\frac{x_1}{x_3}w_1-\frac{x_2}{x_3}w_2=w_1\left(y_1-\frac{x_1}{x_3}\right)+w_2\left(y_2-\frac{x_2}{x_3}\right).$$
Let 
$$ A_1=x_3B_2-x_2B_3, \ A_2=x_1B_3-x_3B_1, \ A_3=B_1-\frac{x_1}{x_3}, \ A_4=B_2-\frac{x_2}{x_3}.$$ 
From the assumption that $\hd(A)>\frac{2}{3}$,  it follows that each $\hd(A_j)$ is also $>\frac{2}{3}$ . By the above calculation, one has 
$$ \{a_1a_2+a_3a_4: a_j \in A_j \} \subset {\mathcal V}(B_1 \times B_2 \times B_3),$$ 
and so the conclusion of Theorem \ref{cartesianproduct} follows by Theorem \ref{tube}. 

\vskip.125in 

\section{Proof of Theorem \ref{crossproduct}} 

\vskip.125in 

As  noted in the introduction, in view of Thm. \ref{witherdogan}, it suffices to prove Thm. \ref{crossproduct}. For clarity and possible future use, we begin the analysis in $\R^d$, specializing to $d=3$ later on. To this end,  define a natural measure on the set of wedge products by the relation 

$$ \int f(z)\,  d\Lambda(z)=\int \dots \int f(x^1 \wedge \dots \wedge x^{d-1}) 
d\mu(x^1) \dots d\mu(x^{d-1}).$$ 

It follows that 

$$ \widehat{\Lambda}(\xi)=\int \dots \int e^{-2 \pi i \xi \cdot *(x^1 \wedge \dots \wedge x^{d-1})} d\mu(x^1) \dots d\mu(x^{d-1})$$ 
and thus 

$$ \int {\left|\widehat{\Lambda}(\xi)\right|}^2 \psi \left( \frac{\xi}{R} \right) d\xi$$ 

\begin{equation} \label{mamalambda}=R^d \int \int \widehat{\psi}(R(x^1 \wedge \dots \wedge x^{d-1}-y^1 \wedge \dots \wedge y^{d-1})) d\mu(x^1) \dots d\mu(x^{d-1}) d\mu(y^1) \dots d\mu(y^{d-1}).\end{equation} 

Since $\widehat{\psi}$ is rapidly decaying, it suffices to estimate 

$$ (\mu \times \dots \times \mu) \big\{(x^1, \dots, x^{d-1}, y^1, \dots, y^{d-1}): |x^1 \wedge \dots \wedge x^{d-1}-y^1 \wedge \dots \wedge y^{d-1}| \leq R^{-1} \big\}.$$

Let $B_r^d(x)$ be the $d$-dimensional ball of radius $r$ centered at $x$; if $x =0$, the center is suppressed. $A^d_{r}$ will denote a $d$-dimensional annulus of inner/outer radii $\frac{r}{2}$ and $r$, not necessarily centered at the origin. For both, $d$ is the ambient dimension if not included in the notation; lower dimensional balls and annuli   in ${\Bbb R}^d$ are denoted using the superscript.
Denote $(d-1)$-tuples of vectors in ${\Bbb R}^d$ by $\vx=(x^1,\dots,x^{d-1})$, and then, as above,  $\hx=*(x^1\wedge\dots\wedge x^{d-1})$, viewed as a vector. 

Let $\mu$ be a Frostman measure on the $s$-dimensional $E\subset A^d_{1}$ \cite{M95}. To  control (\ref{mamalambda}), we want to estimate $R^d\cdot \mu^{2d-2}(F^R)$, where
$$F^R:=\big\{(\vx,\vy): |\hx-\hy|<1/R\big\}.$$
 
Start by decomposing 
$$B^d_1= B_{1/R}(0)\cup\bigcup_{i=0}^{\log_2 R}\, \bigcup_{j=0}^{(R/2^i)^d} \, B_{1/R}(z_j^{(i)}),$$
where $\{z^{(i)}_j\}_j$ is a $(1/R)$-net of points in the dyadic shell $A_{2^{-i-1},2^{-i}}$. Then we can write
\begin{eqnarray*}
 F^R&=& \big\{(\vx,\vy): \hx,\hy\in B_{1/R}\big\}\cup\bigcup_i\bigcup_j\big\{(\vx,\vy):\hx,\hy\in B_{1/R}(z_j^{(i)})\big\} \\
& \subset& (G_0\times G_0) \cup\bigcup_i\bigcup_j G_j^{(i)}\times G_j^{(i)},
\end{eqnarray*} where
$$G_0=\big\{\vx: \hx\in B_{1/R}\big\},\qquad G_j^{(i)}=\big\{\vx: \hx\in B_{1/R}(z_j^{(i)})\big\};$$ 
in terms of  measure, the $G_j^{(i)}$  essentially only depend on $i$  and we refer to them  generically as $G^{(i)}$. With all this, one has
\begin{equation} \label{big3dmama} R^d\mu^{2d-2}(F^R)\lesssim R^d\left(\mu^{d-1}\left(G_0\right)^2 + \sum_{i=0}^{\log_2 R}\left(\frac{R}{2^i}\right)^d\cdot\left(\mu^{d-1}(G^{(i)})\right)^2\right), \end{equation}
where $\mu^k:=\mu\times\dots\times\mu$, and
and we want to estimate the terms on the right hand side.
\bigskip

We now restrict ourselves to three dimensions ($d=3$). 

To estimate the  $G_0$ term, note that $x^1\in E$ is arbitrary, contributing $\mu\lesssim C$. For $x^1$ fixed, $x^2\in $ a $(1/R)\times(1/R)\times 1$ tube, which we cover with $R$ $(1/R)$-balls, giving a $\mu$ contribution $\lesssim R\cdot R^{-s}=R^{1-s}$. Thus, $\mu^2(G_0)\lesssim R^{1-s}$ and hence $\mu^4(G_0\times G_0)\lesssim R^{2-2s}$. 
 
To estimate the $G^{(i)}$ term, start by noting that, with $z_0:=z_j^{(i)}$ fixed, $x^1$ must be in a thin  nhood  of thickness $2^i/R$ in $\R^3$ of an annulus $A^2_1\subset z_0^\perp$; call such a set a washer. Covering this with $(R/2^i)^2$ balls of radius $2^i/R$ gives a $\mu\lesssim (R/2^i)^{2-s}$. For each such $x^1$ fixed, 
 $x^2$ must satisfy two constraints: (i) It has to be in the same washer as $x^1$; and (2) It has to make  an angle of $\sim 2^{-i}$ with $x^1$ and vary in the radial direction by $\le 2^i/R$. There are thus two cases, 
$$ (a) \ 1\le 2^i\le R^\frac12; \ \text{and} \  (b) \ R^\frac12\le 2^i\le R.$$ 
 
In case (a), $2^i/R\le 2^{-i}$, so $x^2$ is confined to a $2^{-i}\times (2^i/R)\times (2^i/R)$-tube, which is covered by  $R/2^{2i}$ $(2^i/R)$-balls, giving 
$$\mu\lesssim (R/2^{2i})(2^i/R)^s=2^{(s-2)i} R^{1-s}.$$ 
Multiplying by the $\mu$ measure in $x^1$ gives 
$$\mu^2(G^{(i)})\lesssim 2^{(2s-4)i}R^{3-2s};$$ squaring this and multiplying  by the prefactor $(R/2^i)^d=R^32^{-3i}$ gives
$$2^{(4s-11)i}R^{9-4s},\quad 1\le 2^i\le R^\frac12,$$ which, since $4s-11<0$ for the $s$  of interest, takes its largest value for  $2^i=1$, corresponding to generic configurations and yielding the estimate $R^{9-4s}$ (before inclusion of  the pre-prefactor of $R^d=R^3$ in (\ref{big3dmama})).
 
In case (b), $2^i/R\ge 2^{-i}$ and so $x^2$ is confined to a $2^{-i}\times 2^{-i}\times (2^i/R)$-tube, which is covered by $2^{2i}/R$ $(2^{-i})$-balls, giving 
$$\mu\le (2^{2i}/R)\cdot 2^{-si}=2^{(2-s)i}R^{-1}.$$ 
Multiplying by the $\mu$ measure in $x^1$, squaring and multiplying by the $R^32^{-3i}$ gives $2^{-3i}R^{5-2s}$ in the range $ R^\frac12\le 2^i\le R$. The largest value is for $2^i=R^\frac12$, namely $R^{\frac72-2s}$, but this is smaller than the $R^{9-4s}$, as is the $G_0$ term.
 
Thus, we estimate (\ref{big3dmama}) by $R^3\cdot R^{9-4s}=R^{12-4s}$. It follows that the expression in (\ref{mamalambda}) is 
$ \lesssim R^{3-(4s-9)}$, which implies that the Hausdorff dimension of the set $\Lambda(E)$
is greater than or equal to $4s-9$. To make use of Thm. \ref{witherdogan}, we need 
$ s+(4s-9)>4$, which holds if $s>\frac{13}{5}$. This completes the proof of Thm. \ref{crossproduct} and thus of Thm. \ref{maingeneral}. 

\vskip.25in 

\section{Applications to discrete geometry} 

\vskip.125in 

\begin{definition} \cite{IRU10} Let $P$ be a set of $n$ points contained in ${[0,1]}^2$. Define the measure
\begin{equation} \label{pizdatayamera} d \mu^s_P(x)=n^{-1} \cdot n^{\frac{d}{s}} \cdot \sum_{p \in P} \chi_{B^{\, p}_{n^{-\frac{1}{s}}}}(x)dx, \end{equation} where $\chi_{B^{\, p}_{n^{-\frac{1}{s}}}}(x)$ is the characteristic function of the ball of radius $n^{-\frac{1}{s}}$ centered at $p$. 
One says that $P$ is \emph{$s$-adaptable} if 

$$ I_s(\mu_P)=\int \int {|x-y|}^{-s} d\mu^s_P(x) d\mu^s_P(y)<\infty.$$ 

\end{definition} 

\vskip.125in

This is equivalent to the statement 

\begin{equation} \label{discreteenergy} n^{-2} \sum_{p \not=p' \in P} {|p-p'|}^{-s} \lesssim 1.\end{equation} 

\vskip.125in 

To understand this condition in  clearer geometric terms, suppose that $P$ comes from a \linebreak$1$-separated set $A$, scaled down by its diameter. Then the condition (\ref{discreteenergy}) takes the form 

\begin{equation} \label{discreteenergylarge} n^{-2} \sum_{a \not=a' \in A} {|a-a'|}^{-s} \lesssim {(diameter(A))}^{-s}. \end{equation} 

Thus, $P$ is $s$-adaptable if it is a scaled $1$-separated set where the expected value of the distance between two points raised to the power $-s$ is comparable to the value of the diameter raised to the power of $-s$. This basically means that clustering within $P$ is not allowed to be too severe. 

In more technical terms, $s$-adaptability means that a discrete point set $P$ can be thickened into a set which is uniformly $s$-dimensional in the sense that its energy integral of order $s$ is finite. Unfortunately, it is shown in \cite{IRU10} that there exist finite point sets which are not $s$-adaptable for certain ranges of the parameter $s$. The point is that the notion of Hausdorff dimension is much more subtle than the simple ``size" estimate. However, many natural classes of sets are $s$-adaptable. For example, homogeneous sets studied by Solymosi and Vu \cite{SV04} and others are $s$-adaptable for all $0<s<d$. See also \cite{IJL09} where $s$-adaptability of homogeneous sets is used to extract discrete incidence theorems from Fourier-type bounds.

The following argument is a variant of the conversion mechanism developed in \cite{IL05,HI05}. Suppose that  one knew that ${\mathcal L}^1({\mathcal V}(E))>0$ whenever $\hde>s_\V$. Let $P$ be an $s$-adaptable set with $s \in (s_0,d]$, and $E$ denote the support of $d\mu^s_P$ above. It follows that 
$$ 1 \lesssim {\mathcal L}^1({\mathcal V}(E)) \lesssim n^{-\frac{1}{s}} \cdot \# {\mathcal V}(P),$$ 
and one can conclude that 
$$ \# {\mathcal V}(E) \gtrsim n^{\frac{1}{s}} \gtrapprox n^{\frac{1}{s_0}},$$ which yields the following theorem. 

\begin{theorem} \label{discretevolume} Let $P$ be a $\frac{13}{5}$-adaptable subset of ${\Bbb R}^d$ of size $n$. Let 
${\mathcal V}_{\delta}(P)$ denote the number of distinct $\delta$-separated volumes determined by $P$. Then 
$$ \# {\mathcal V}_{n^{-\frac{5}{13}}}(P) \gtrapprox n^{\frac{5}{13}}.$$ 
\end{theorem} 

It is important to note that the significance of this result lies in the requirement that the volumes are $n^{-\frac{5}{13}}$-separated.
In the absence of this feature, the exponent is  inferior to the sharp result obtained by Dumitrescu and Toth \cite{DT11}, who proved that for any points set $P \subset {\Bbb R}^3$ of size $n$, $\# {\mathcal V}(P) \gtrapprox n$. 

Even if one were to improve  Thm. \ref{maingeneral} to the conjecturally sharp $\hde>2$, 
the resulting modification to Thm. \ref{discretevolume} would be that $\# {\mathcal V}_{n^{-\frac{1}{2}}}(P) \gtrapprox n^{\frac{1}{2}}$, still weaker in terms of the exponent than \cite{DT11}. The analytic method used here yields a conclusion about $n^{-\frac{1}{s}}$-separated volumes; we do not know whether it is possible to obtain the sharp exponent  using our methods. 


\end{document}